%genpark.tex: a Plain TeX file by AJ Bu and Doron Zeilberger
%A  Short  Proof that the number of (a,b)-parking functions of length n is a(a+bn)^(n-1) 

%begin macros

\baselineskip=14pt
\parskip=10pt
\def\halmos{\hbox{\vrule height0.15cm width0.01cm\vbox{\hrule height
  0.01cm width0.2cm \vskip0.15cm \hrule height 0.01cm width0.2cm}\vrule
  height0.15cm width 0.01cm}}
\font\eightrm=cmr8 

\magnification=\magstephalf

\def\B{{\cal B}}

\def\P{{\cal P}}

\def\1{{\overline{1}}}
\def\2{{\overline{2}}}
\parindent=0pt
\overfullrule=0in

\def\frac#1#2{{#1 \over #2}}
%\headline={\rm  \ifodd\pageno  \RightHead  \else  \LeftHead  \fi}
%\def\RightHead{\centerline{
%Title
%}}
%\def\LeftHead{ \centerline{Doron Zeilberger}}
%end macros
\centerline
{\bf 
A  Short  Proof that the number of $(a,b)$-parking functions of length $n$ is $a(a+bn)^{n-1}$
}
\bigskip
\centerline
{\it AJ BU and Doron ZEILBERGER}
\bigskip

{\bf Abstract}: We give a very short proof of the fact that the number of $(a,b)$-parking functions of length $n$ equals $a(a+bn)^{n-1}$. This was first proved
in 2003 by Kung and Yan, via a very long and torturous route, as a corollary of a more general result.

Recall [KoW] that a parking function of length $n$ is a list of positive integers $x=(x_1, \dots, x_n)$ whose  non-decreasing rearrangement, $x'_1 \leq \dots \leq x'_n$
satisfies $x'_i \leq i, (1 \leq i \leq n)$. Their number is famously $(n+1)^{n-1}$. 

More generally, for any non-negative integers, $a$ and $b$, an $(a,b)$-parking function of length $n$ is 
a list of positive integers $(x_1, \dots, x_n)$ whose  non-decreasing rearrangement satisfies:
$x'_i \leq a+b(i-1), (1 \leq i \leq n)$. Note that $(1,1)$-parking functions are the usual ones.

In [KuY], after $13$ pages of heavy-going math, as a corollary of a much more general result (Cor. 5.5 there), the authors proved the following fact.

{\bf Fact}: The number of $(a,b)$-parking functions of length $n$ is $a(a+bn)^{n-1}$.

{\bf Short Proof}: Let $\P(n,a,b)$ be the set of $(a,b)$-parking functions of length $n$, and let $p(n,a,b):=|\P(n,a,b)|$. For $0 \leq r \leq n$, the subset consisting of those with exactly $r$ $1$s is in bijection with
$\B_{n,r}\times \P(n-r,a+br-1,b)$, where $\B_{n,r}$ is the set of $r$-element subsets of $\{1, \dots, n\}$.
It is given by: $ x \rightarrow (S,y)$, where $S$ is the subset of locations of $x$ where the $1$s reside, and $y$ is
obtained from $x$ by deleting the $1$s, and then subtracting one from every remaining entry, and since $x'_{r+1} -1 \leq a\,+\,b\,r\,-\,1$, the claim follows. Hence:
$$
p(n,a,b)= \sum_{r=0}^{n} {{n} \choose {r}} p(n-r,a+br-1,b) \quad ,
$$
that uniquely determines $p(n,a,b)$, subject to the initial conditions $p(0,a,b)=1$ and $p(n,0,b)=0$.

But $q(n,a,b):=a(a+bn)^{n-1}$ satisfies the very same recurrence, where $p$ is replaced by $q$, namely:
$$
a(a+bn)^{n-1}= \sum_{r=0}^{n} {{n} \choose {r}} (a+br-1)(a+bn-1)^{n-r-1} \quad ,
$$
(Check!\footnote{${}^1$}
{\eightrm
Write $a\,+\,b\,r\,-\,1$ as $(a\,+ \, b\,n \,- \,1)\,-\,b \,(n\,-\,r)$ , use $(n-r){{n} \choose {r}}=n{{n-1} \choose {r}}$ and then invoke the binomial theorem twice.}), 
and hence $p(n,a,b)=q(n,a,b)$ follows by induction on $n$. \halmos

\vfill\eject

{\bf Comment:}

We were unaware of this result, and found it {\it ab initio}, by using {\it experimental mathematics}, via the Maple package

{\tt https://sites.math.rutgers.edu/\~{}zeilberg/tokhniot/GenPark.txt }. We thank Lucy Martinez for telling us about [KuY].

{\bf Added Dec. 20, 2024:} It turns out that this is {\bf not} the first short proof. Richard Stanely emailed us with the following:

In connection with your recent arXiv paper with AJ Bu, you might be
interested in the paragraph preceding Theorem 1.2 of the paper 

Richard P. Stanley and Yinghui Wang, {\it Some aspects of $(r,k)$-parking functions}, J. Combinatorial Theory (A) {\bf 159 }(2018), 54-78,
{\tt https://math.mit.edu/\~{}rstan/papers/pf.pdf} \quad .

It also turns out (see reference [7] in the above paper) that this result is much older than [KuY], it goes back to 1969  and is due to G. P. Steck.

{\bf References}

[KoW] Alan G. Konheim and Benjamin Weiss, {\it An occupancy discipline and applications},
SIAM J. Applied Math. {\bf 14} (1966), 1266-1274. [Available from JSTOR.] 

[KuY] Joseph P.S. Kung and Catherine Yan, {\it Gon{\u c}arov polynomials and parking functions}, J. of Combinatorial Theory Ser. A, {\bf 102} (2003), 16-37.

\bigskip
\hrule
\bigskip
AJ Bu and Doron Zeilberger, Department of Mathematics, Rutgers University (New Brunswick), Hill Center-Busch Campus, 110 Frelinghuysen
Rd., Piscataway, NJ 08854-8019, USA. \hfill\break
Email: {\tt   ajbu131@gmail.com} \quad, \quad {\tt DoronZeil@gmail.com}   \quad .

First Written: {\bf Dec. 15, 2024}.

This version: {\bf Dec. 20, 2024}.
\end